\documentclass[a4paper]{article}
\usepackage{a4wide}
\usepackage{url,hyperref,xspace,amssymb,graphicx,hyperref}
\usepackage[comma]{natbib}
\newcommand{\ie}{{\em i.e.\/}\xspace}
\newcommand{\eg}{{\em e.g.\/}\xspace}

\newcommand{\HIDE}[1]{ }
\newcommand{\COMMENT}[1]{ }

\newcommand{\pc}{{\sc PC}\xspace}

\renewcommand{\dag}{{\cal D}}

\newcommand{\eqref}[1]{\mbox{(\ref{eq:#1})}}

\newcommand{\secref}[1]{\mbox{\S$\,$\ref{sec:#1}}}

\newcommand{\figref}[1]{\mbox{Figure~\ref{fig:#1}}}

\newcommand{\itref}[1]{\mbox{\ref{it:#1}}}

\newcommand{\exref}[1]{\mbox{Example~\ref{ex:#1}}}

\newcommand{\Secref}[1]{\mbox{Section~\ref{sec:#1}}}

\newcommand{\cip}{\mbox{$\perp\!\!\!\perp$}}
\newcommand{\indo}[2]{\mbox{$#1 \,\cip\, #2$}}
\newcommand{\ind}[3]{\mbox{$#1 \, \cip\, #2 \mid #3$}}

\newtheorem{expl}{Example}%[section]
\newtheorem{definer}{Definition}%[section]
\newtheorem{algor}{Algorithm}%[section]
\newtheorem{rem*}{Remark}%[section]
\newcommand{\halm}{\hspace*{\fill} $\Box$\par}
\newenvironment{ex}{\begin{expl}\rm}{\halm\end{expl}}

\newcommand{\pr}[1]{{\rm prob}(#1)} 
\renewcommand{\pc}{\mbox{\rm PC}\xspace}
\renewcommand{\pr}{\mbox{\rm P}\xspace}
\newcommand{\bY}{\mbox{\boldmath$Y$}}

\title{What Is a Causal Graph?}

\author{Philip Dawid\footnote
{University of Cambridge}}
\date{}
\begin{document}
\maketitle

%%%%%%%%%%%%%%%%%%%%%%%%%%%%%%%%%%%%%%%%%%
\abstract{This article surveys the variety of ways in which a directed
acyclic graph (DAG) can be used to represent a problem of
probabilistic causality.  For each of these we describe the relevant
formal or informal semantics governing that representation.  It is
suggested that the cleanest such representation is that embodied in an
augmented DAG, which contains nodes for non-stochastic intervention
indicators in addition to the usual nodes for domain variables.\\

}
\section{Introduction}
\begin{figure}[htbp]
  % \hspace{1cm} {\epsfbox{fig1.eps}[width=3cm]}
\begin{center}
  \resizebox{3in}{!}{\includegraphics{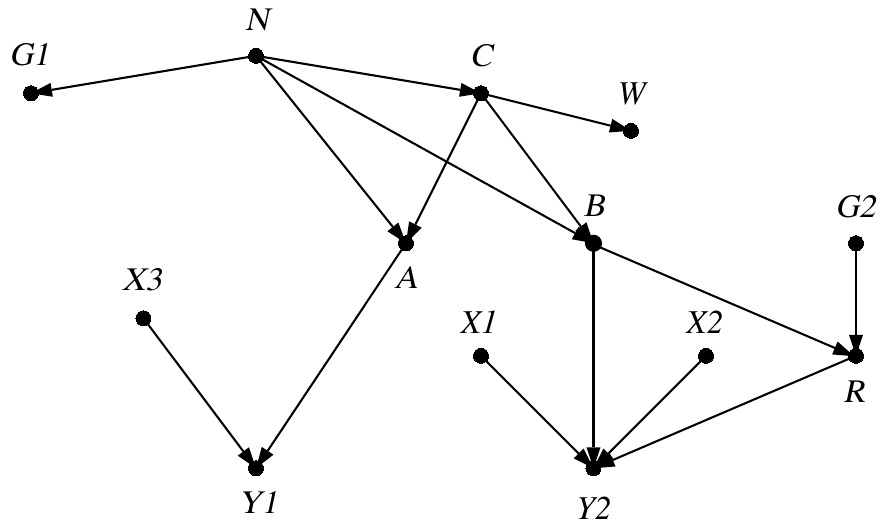}}
\caption{Directed acyclic graph $\dag$}
\label{fig:dag}
\end{center}
\end{figure}

Graphical representations of probabilistic
\citep{pearl:pesbook,sll:book,mybook} and causal \citep{pearl:book}
problems are ubiquitous.  Such a graph has nodes representing relevant
variables in the system, which we term {\em domain variables\/}, and
arcs between some of the nodes.  The most commonly used type of graph
for these purposes, to which we will confine attention in this
article\footnote{However other types of graph also have valuable
  applications}, is a Directed Acyclic Graph (DAG), in which the arcs
are arrows, and it is not possible to return to one's starting point
by following the arrows.  An illustration of such a DAG \citep{apd/iwe}
is given in \figref{dag}.

Now there is no necessary relationship between a geometric object,
such as a graph, and a probabilistic or causal model---they inhabit
totally different mathematical universes.  Any such relationship must
therefore be specified externally, and then constitutes a way of
interpreting the graph as saying something about the problem at hand.
It is important to distinguish between the {\em syntax\/} of a graph,
\ie, its internal, purely geometric properties, and its {\em
  semantics\/}, describing its intended interpretation.

In this article we consider a variety of ways---some formal, some less
so---in which DAGs have been and can be interpreted.  In particular,
we emphasise the importance of a clear understanding of what is the
intended interpretation in any specific case.  We caution against the
temptation to interpret purely syntactical elements, such as arrows,
semantically (the sin of ``reification''), or to slide unthinkingly
from one interpretation of the graph to another.

The material in this article is largely a survey of previously
published material
\citep{apd:infdiags,apd:hsss,apd:beware,apd:pearl,apd:annrev,apd:found,apd:wag,dd:ss,sgg/apd:ett,apd/mm/rm,apd/mm:eoccoe,pc/apd:eci,dhm},
which should be consulted for additional detail and discussion.  In
particular we do not here showcase the wealth of important
applications of the methods discussed.

\subsection*{Outline}
\label{sec:outline}
In \secref{probdag} we describe how a DAG can be used to represent
properties of probabilistic independence between variables, using a
precise semantics based on the method of moralisation.  An alternative
representation, involving functional rather than probabilistic
dependencies, is also presented.  \Secref{inform} discusses some
informal and imprecise ways in which a DAG might be considered as
representing causal relations.  In preparation for a more formal
approach, \secref{int} describes an understanding of causality as a
reaction to an external intervention in a system, and presents an
associated language, based on an extension of the calculus of
conditional independence, for expressing and manipulating causal
concepts.  Turning back to a DAG model, in \secref{augdag} we
introduce a precise semantics for its causal interpretation, again
based on moralisation, this time used to express extended condtional
independence properties.  To this end we introduce augmented DAGs,
which contain nodes for non-stochastic intervention indicators, as
well as for stochastic domain variables.  In \secref{pearl} we
describe the causal semantics of a ``Pearlian'' \citep{pearl:book} DAG,
and show how these can be explicitly represented by an augmented DAG,
where again a version involving functional relationships---the
Structural Causal Model (SCM) ---is available.

\Secref{coe} considers a different type of causal problem, that is not
about the response of a system to an intervention, but rather aims to
assign, to an observed putative cause, responsibility for the
emergence of an observed outcome.  This requires new understandings
and semantics that can not be represented by a probabilistic DAG, but
can be by using a structural causal model.  However this is
problematic, since there can be distinct SCMs that are observationally
equivalent but lead to different answers.

\section{Probabilistic DAG}
\label{sec:probdag}

The most basic way of interpreting a DAG is as representing
qualitative probabilistic relations of conditional independence
between its variables. Such a situation occupies the lowest rung of
the ``ladder of causaton'' \citep{pearl:why}.  The semantics governing
such a representation, while precise, are not totally straightforward,
being described by either of two logically equivalent criteria, known
as ``$d$-separation'' \citep{pearl:blocking,verma} and ``moralisation''
\citep{ldll}.  Here we describe moralisation, using the specific DAG
$\dag$ of \figref{dag} for illustration.

Suppose we ask: Is the specific conditional independence property
$$\ind {(B,R)}{(G_1,Y_1)}{(A,N)}$$
(read as ``the pair $(B,R)$ is independent of $(G_1,Y_1)$, conditional
on $(A,N)$''---see \cite{apd:CIST}) represented by the graph?  To
address this query we proceed as follows:
\begin{description}
\item[1.  Ancestral graph]

  We form a new DAG $\dag'$ by removing from $\dag$ every node that is
  not mentioned in the query and is not an ancestor of a mentioned
  node\footnote{\ie, there is no directed path from it to a mentioned
    node}, as well as any arrow involving a removed node.
  \figref{anc} shows the result of applying this operation to
  \figref{dag}.

  \begin{figure}
  % \hspace{1cm} {\epsfbox{fig2}}
  \begin{center}
    \resizebox{3in}{!}{\includegraphics{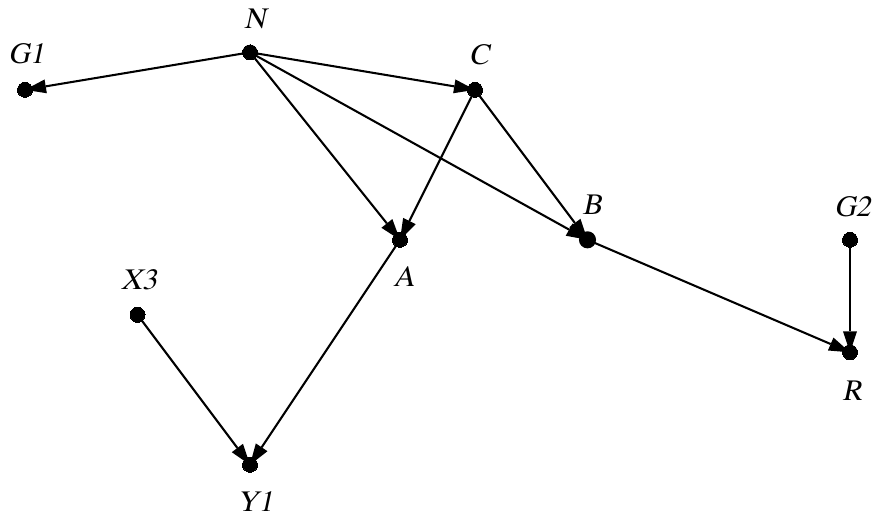}}
  \end{center}
\caption{Ancestral subgraph ${\cal D}'$}
\label{fig:anc}
\end{figure}

\item[2. Moralisation] If in $\dag'$ there are two nodes that have a
  common child (\ie, each has an arrow directed from it to the child
  node) but they are not themselves joined by an arrow---a
  configuration termed an ``immorality''---then an undirected edge is
  inserted between them.  Then every remaining arrowhead is removed,
  yielding an undirected graph ${\cal G}'$.  In our example this
  yields \figref{mor}.

  \begin{figure}%\hspace{1cm}
    % {\epsfbox{fig3.eps}}
    \begin{center}
      \resizebox{3in}{!}{\includegraphics{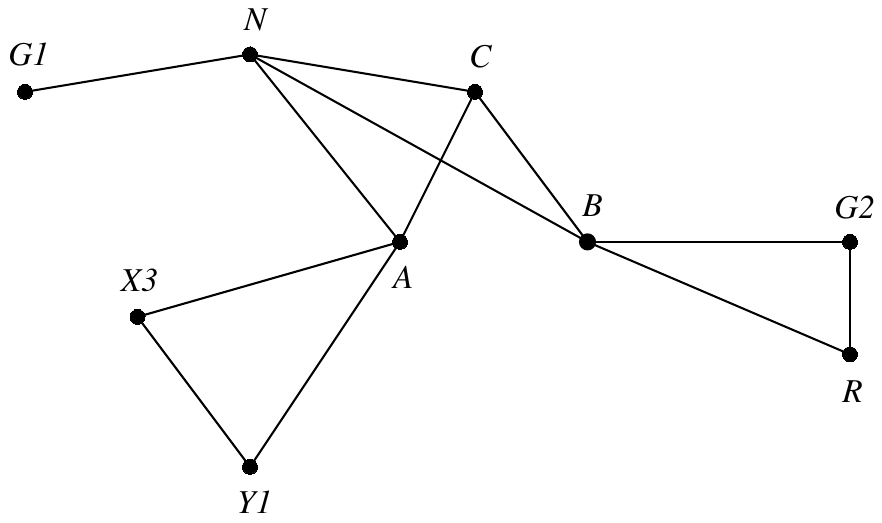}}
    \end{center}
    \caption{Moralised ancestral subgraph ${\cal G}'$}
\label{fig:mor}
\end{figure}

\item[3.  Separation] Finally, in ${\cal G}'$ we look for a continuous
  path connecting an element of the first set in our query (here
  $(B,R)$) to an element of the second set (here $(G_1,Y_1)$) that
  does not contain an element of the third set (here $(A,N)$).  If
  there is no such path, we conclude that the queried conditional
  independence property is represented in the original DAG.  Since
  this is the case in our example, the property
  $\ind {(B,R)}{(G_1,Y_1)}{(A,N)}$ is indeed represented in
  \figref{dag}.
\end{description}

\subsection{Instrumental variable}
\label{sec:instr}

\begin{figure}[htbp]
  \centering
  \includegraphics [width=.4\textwidth,clip] {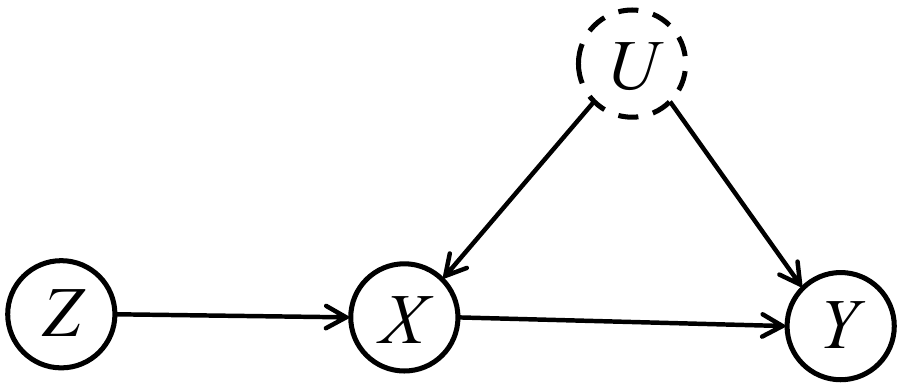}
  \caption{Instrumental variable}
  \label{fig:instr}
\end{figure}

The DAG of \figref{instr} can be  used to represent a problem in which we have an ``instrumental variable'', with the nodes interpreted as follows:

\begin{description}
\item[$X$] Exposure
\item[$Y$] Outcome
\item[$Z$] Instrument
\item[$U$] Unobserved confounding variables\footnote{The dotted
    outline on node $U$ serves as a reminder that $U$ is unobsderved,
    but is otherwise of no consequence.}
\end{description}

Such a problem invites a causal interpretation, which we will take up
in \secref{augdag} below.  But if we interpret \figref{instr} purely
as a probabilistic DAG, we see that it represents exactly the
following conditional independencies:
\begin{eqnarray}
  \label{eq:instrp1}
  U &\cip& Z\\
  \label{eq:instrp2}
  Y &\cip& Z \mid (U,X)
\end{eqnarray}

\subsection{Markov equivalence} Given a collection of conditional
independence relations for a set of variables, there may be $0$, $1$,
or several DAGs that represent just these relations.  Two DAGs are
termed {\em Markov equivalent\/} when they represent the identical
conditional independencies.  It can be shown \citep{morten:markov,
  verma/pearl:90} that this will be the case if they have the same
skeleton (\ie, the undirected graph obtained by deleting the
arrowheads) and the same immoralities.

\begin{ex}
  \label{ex:DAG-markeqsimp}
  Consider the following DAGs on three nodes:
  \begin{enumerate}
\item 
      \label{it:DAG-markeqsimp1}
      $A \to B \to C$
\item 
      \label{it:DAG-markeqsimp2}
      $A \leftarrow B \leftarrow C$
\item 
      \label{it:DAG-markeqsimp3}
      $A \leftarrow B \to C$
\item 
      \label{it:DAG-markeqsimp4}
      $A \to B \leftarrow C$.
  \end{enumerate}
  These all have the same skeleton.  However, whereas DAGs
  \itref{DAG-markeqsimp1}, \itref{DAG-markeqsimp2} and
  \itref{DAG-markeqsimp3} have no immoralities,
  \itref{DAG-markeqsimp4} has one immorality.  Consequently,
  \itref{DAG-markeqsimp1}, \itref{DAG-markeqsimp2} and
  \itref{DAG-markeqsimp3} are all Markov equivalent, but
  \itref{DAG-markeqsimp4} is not Markov equivalent to these.  Indeed,
  \itref{DAG-markeqsimp1}, \itref{DAG-markeqsimp2} and
  \itref{DAG-markeqsimp3} all represent the conditional independence
  property $\ind A C B$, whereas \itref{DAG-markeqsimp4} represents the
  marginal independence property $\indo A C$.
\end{ex}

\subsection{Bayesian network}
\label{sec:bayes}
The purely qualitative graphical structure of a probabilistic DAG can
be elaborated with quantitative information.  With each node in the
DAG we associate a specified conditional distribution for its
variable, given any values for its parent variables.  There is a
one-one correspondence between a collection of all such parent-child
distributions and a joint distribution for all the variables that
satisfies all the conditional independencies represented by the graph.
The graphical structure also supports elegant algorithms for computing
marginal and conditional distributions \citep{mybook}.

\subsection{Structural probabilistic model}
\label{sec:struct}
Suppose we are given the conditional distribution $p(Y \mid X)$.  It
is then possible to construct (albeit non-uniquely) a fictitious
``error variable'' $E$, having a suitable distribution $P$, and a
suitable deterministic function $f$ of $(X,E)$, such that the
distribution of $f(x, E)$ is just $p(Y \mid X=x)$.  For example, using
the probability integral transform, if $F_x$ is the cumulative
distribution function of $p(Y \mid X=x)$, we can take
$f(x,e) = F_x^{-1}(e)$\footnote{using a suitable definition of the
  inverse function $F_x^{-1}$ when $F_x$ is not continuous} and $E$
uniform on $[0,1]$.  However such a representation is highly
non-unique.  Indeed, we could alternatively take
$f(x,e) = F_x^{-1}(e_x)$, where $e$ is now a vector $(e_x)$, and the
multivariate distribution of $E$ is an arbitrarily dependent copula,
such that each entry $(E_x)$ is uniform on $[0,1]$.

Given any such construction, for purely distributional purposes we can
work with the functional equation $Y=f(X,E)$, with $E\sim P$
independently of $X$.  This operation can be extended to apply to a
complete Bayesian network, by associating with each domain variable
$V$ a new error variable $E_V$, with suitable distribution $P_v$, all
these being independent, and with $V$ being modelled as a suitable
determistic function of its graph parent domain variables and $E_V$.
We term such a deterministic model a Structural Probabilistic Model
(SPM), by analogy with the Structural Causal Model (SCM)
\citep{pearl:why}---see \secref{structcaus} below.  In an SPM, all
stochasticity is confined to the error variables.

\begin{figure}[htbp]
  \centering
  \includegraphics [width=.4\textwidth,clip] {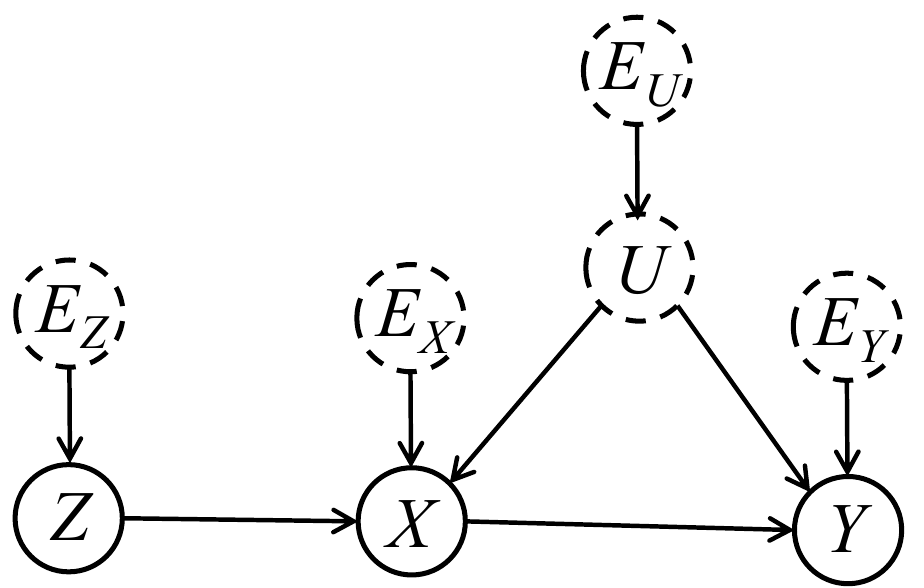}
  \caption{SPM representation of instrumental variable model}
  \label{fig:instre}
\end{figure}

An SPM can be represented graphically by introducing new nodes for the
error variables, as illustrated in \figref{instre} for the case of
\figref{instr}, it being understood that the error variables are
modelled as random, but all other parent-child relations are
deterministic.

It is easy to check that the conditional independencies represented
between the domain variables, as revealed by the moralisation
criterion, are identical, both in the original graph and in an SPM
extension of it.  So, when each graph is endowed with its own
parent-child relationships (stochastic in the case of the original
graph, deterministic for the structural graph), both graphs describe
the identical joint distribution for the domain variables.  For purely
probabilistic purposes, nothing is gained by introducing error
variables to ``mop up'' the stochastic dependencies of a simple DAG
model.

\subsection{Reification?}
\label{sec:reification}
It is important to keep in mind that a probabilistic DAG is nothing
but a very indirect way of representing a collection of conditional
independence properties.  In particular, the arrows in the DAG have no
intrinsic meaning, being there only to support the moralisation
procedure.  It is indeed somewhat odd that the property of conditional
independence, which is essentially a symmetric one, can be represented
at all by means of arrows, which have built-in directionality.  The
example of Markov equivalence between \itref{DAG-markeqsimp1},
\itref{DAG-markeqsimp2} and \itref{DAG-markeqsimp3} in
\exref{DAG-markeqsimp} shows that the specific direction of an arrow
in a DAG should not be taken as meaningful itself.  Rather, an arrow
has a transient status rather like that of a construction line in an
architect's plan, or of a contour line on a map: instrumentally
helpful, but not representing anything visible in the house or on the
ground.

Nevertheless, on looking at a DAG, it is hard to avoid the temptation
to imbue its arrows with a direct meaning in relation to the system
studied.  This is the philosophical error of {\em reification\/},
which confuses the map with the territory \citep{korzybski}, wrongly
interpreting a purely instrumental property of a representation as if
it had a direct counterpart in the external system.  In the case of a
probabilistic DAG model, this frequently leads to endowing it with an
unjustified causal interpretation, where the presence of an arrow
$A\rightarrow B$ is considered to represent the presence of a causal
influence of $A$ on $B$.  Such a confusion underlies much of the
enterprise of ``causal discovery'', where observational data are
analysed to uncover their underlying probabilistic conditional
indepencies, these are represented by a probabilistic DAG, and that
DAG is then reinterpreted as representing causal relations.

This is not to say that a DAG can not be interpreted causally---but to
do so will require a fresh start, with new semantics.

\section{Informal causal semantics}
\label{sec:inform}
Common causal interpretations of a DAG involve statements such as:
\begin{itemize}
\item an {\bf arrow} represents a {\tt direct cause}
\item a {\bf directed path} represents a {\tt causal pathway}
\end{itemize}
Or, as described for example by \cite{hernan:06}:
\begin{quotation}
  ``A {\em causal DAG\/} $\dag$ is a DAG in which:
  \begin{enumerate}
  \item \label{it:hr1} the {\bf lack of an arrow} from $V_j$ to $V_m$ can be
    interpreted as the absence of a {\tt direct causal effect} of
    $V_j$ on $V_m$ {\tt (relative} {to} {the} {other} {variables} {on}
    {the} {graph)}
  \item \label{it:hr2} all {\tt common causes}, even if
    unmeasured, of any pair of variables on the graph are themselves
    {\bf on the graph}.
  \end{enumerate}
\end{quotation} 

Thus \figref{instr} might be interpreted as saying:
\begin{itemize}
\item $U$ is a {\tt common cause} of $X$ and $Y$
\item $Z$ {\tt affects} the outcome $Y$ {\tt only through $X$}
\item $Z$ does not share {\tt common causes} with $Y$
\item $Z$ has a {\tt direct causal effect} on $X$
\end{itemize}
\vspace{1ex} In the above, we have marked syntactical terms, relating
to the graph itself (the ``map''), in {\bf bold-face}, and terms
involving external concepts (the ``territory''), in {\tt teletype}.

\subsection{Probabilistic Causality}
\label{sec:prob-caus}
Or, one can start by ignoring the DAG, and try to relate causal
concepts directly to conditional independence.  Such an approach
underlies the enterprise of {\em causal discovery\/}, where discovered
conditional independencies are taken to signify causal relations.
Common assumptions made here are:\footnote{Even when all terms
  involved are fully understood, one might question just why these
  assumptions should be regarded as appropriate}
\begin{description}
\item[Weak Causal Markov assumption:]
If $X$ and $Y$ have no {\tt common cause} (including each other), they
are {\bf probabilistically independent}
\item[Causal Markov assumption:]
A variable is {\bf probabilistically independent} of its {\tt
  non-effects}, given its {\tt direct causes}
\end{description}
Conbining such assumptions with the formal semantics by which a DAG
represents conditional independence, we obtain a link between the DAG
and causal concepts.

\subsection{A problem}
How could we check if a DAG, endowed with a causal interpretation
along the above lines, is a {\em correct representation\/} of the
external problem it is intended to model?  To even begin to do so, we
would already need to have an independent understanding of the
external concepts, as marked in bold-face---which thus can not make
reference to the graph itself.  But these informal causal concepts are
woolly and hard to pin down.  I am aware of almost no work that even
attempts to do so---a recent exception is \cite{kayvan:ax}.  But
without this, using such informal causal semantics risks generating
confusion rather than clarification.

To avoid confusion, we need to develop a more formal causal semantics.

\section{Interventional causality and extended conditional
  independence}
\label{sec:int}
Our approach to this begins by introducing an explicit, non-graphical,
understanding of causality, expressed in terms of the probabilistic
response of a system to an (actual or proposed) intervention.  A {\tt
  causal effect} of $A$ on $B$ exists if the distribution of $B$,
after an {\em external intervention\/} sets $A$ to some value $a$,
varies with $a$.  Formally we introduce a non-stochastic variable
$F_A$, having the same set of values as $A$, such that $F_A=a$
describes the regime in which $A$ is set to value $a$ by an external
intervention.\footnote{We shall here only consider ``surgical
  interventions'', such that $F_A = a \Rightarrow A=a$.} Then $A$ has
{\tt no causal effect} on $B$ just when the distribution of $B$, given
$F_A=a$, does not depend on $a$.

If $F_A$ were a stochastic variable, this would just be the usual
property of independence of $B$ from $F_A$, notated as
\begin{equation}
  \label{eq:causeff}
  \indo B {F_A}.
\end{equation}
Now not only does this understanding of independence remain
intuitively meaningful in the current case that $F_A$ is a
non-stochastic variable, but the formal mathematical properties of
independence and conditional independence can be rigorously extended
to such cases \citep{apd:CIST,apd:ciso,nayia/apd:eci}: we term this
{\em extended conditional independence (ECI)\/}.  The standard axioms
of conditional independence apply, with essentially no
changes\footnote{In ECI we interpret $\ind A B C$ as the property that
  the conditional distribution of $A$, given $B=b$ and $C=c$, depends
  only on $c$.  Here we can allow $B$ and $C$ to include
  non-stochastic variables; however, $A$ must be fully stochastic.},
to ECI.

Since much of the enterprise of statistical causality involves drawing
causal conclusions from observational data, we further endow $F_A$
with an additional state $\emptyset$, read as ``idle'':
$F_A = \emptyset$ denotes the regime in which no intervention is
applied to $A$, but it is allowed to develop ``naturally'', in the
observational setting.  The distinction between $A$ and $F_A$ is that
the state of $A$ describes {\em what\/} value was taken, while the
state of $F_A$ describes {\em how\/} that value was taken.

We regard the main thrust of causal
statistical inference as aiming to deduce properties of a hypothetical
interventional regime, such as the distribution of $B$ given $F_A=a$,
from observational data, obtained under the idle regime
$F_A = \emptyset$.  But since there is no logically necessary
connexion between the distributions under different regimes, suitable
assumptions will need to be made---and justified---to introduce such
connexions.

The simplest such assumption---which, however will only very rarely be
justifiable---is that when we consider the distribution of $B$ given
$A$, we need to know {\em what\/} was the value $a$ that $A$ took, but
not {\em how\/} it came to take that value (\ie, not whether this was
in the observational regime $F_A=\emptyset$, or in the interventional
regimes, $F_A = a$), the conditional distribution of $B$ being the
same in both cases.  That is,
\begin{equation}
  \label{eq:ig}
  p(B=b \mid A=a, F_A = \emptyset) = p(B=b \mid A=a, F_A = a).
\end{equation}
This is the property of {\em ignorability}.  When this strong property
holds, we can directly take the observed distribution of $B$ given
$A=a$ for the desired interventional distribution of $B$ given
$F_A=a$: the distribution of $B$ given $A$ is a ``modular component'',
transferable between different regimes.

We note that \eqref{ig} can be expressed succinctly using ECI notation,
as
\begin{equation}
  \label{eq:igci}
  \ind B {F_A} A
\end{equation}
As exemplified by \eqref{causeff} and \eqref{igci}, ECI supplies a
simple and powerful language for expressing and manipulating causal
concepts.

\section{Augmented DAG}
\label{sec:augdag}
We can now introduce a formal graphical causal semantics, so moving
onto the second rung of the ladder of causation.  This is based on
ECI, represented by an {\em augmented DAG\/}, which is just like a
regular DAG, except that some of its nodes may represent
non-stochastic variables, such as regime indicators.\footnote{We may
  indicate a non-stochastic variable by a square node, and a
  stochastic variable by a round node.  However this distinction does
  not affect how we use the DAG.}  Just as certain collections of
purely probabilistic conditional independence properties can usefully
and precisely be represented (by moralisation) by means of a regular
DAG, so we may be able to construct an augmented DAG to represent (by
exactly the same moralisation criterion) causal properties of interest
expressed in terms of ECI.

Consider for example the simple augmented DAG
$$F_A\rightarrow A \rightarrow B$$
This represents (by moralisation, as always) the ECI property of
\eqref{igci}, and so is a graphical representation of the ignorability
assumption.  Note that it is the whole structure that, with ECI,
imparts causal meaning to the arrow from $A$ to $B$---the direction of
that arrow would not otherwise be meaningful in itself.

\begin{ex}
  \label{ex:augabc}
  Consider the following augmented DAGS, modifications of the first
  three DAGs of \exref{DAG-markeqsimp} to allow for an intervention on
  $A$:
\begin{enumerate}
\item
  \label{it:augsimp1}
  $F_A \to A \to B \to C$
\item 
  \label{it:augsimp2}
  $F_A \to A \leftarrow B \leftarrow C$
\item 
  \label{it:augsimp3}
  $F_A \to A \leftarrow B \to C$
% \item 
%   \label{it:augsimp4}
%   $F_A \to A \to B \leftarrow C$.
\end{enumerate}
We saw in \exref{DAG-markeqsimp} that, without the node $F_A$, these
DAGs are all Markov equivalent.  With $F_A$ included, we see that DAGs
\itref{augsimp2} and \itref{augsimp3} are still Markov equivalent,
since they have the same skeleton and the same unique immorality
$F_A \to A \leftarrow B$; but they are no longer Markov equivalent to
DAG \itref{augsimp1}, which does not contain that immorality.  All
three DAGs represent the ECI $\ind A C {(F_A, B)}$, which says that,
in any regime\footnote{In fact this only has bite for the idle regime,
  since, in an interventional regime $F_A = a$, $A$ has a degenerate
  distribution at $a$, so that the conditional independence is
  trivially satisfied.}, $A$ is independent of $C$ given $B$; however,
while DAGs \itref{augsimp2} and \itref{augsimp3} both represent the
ECI $\indo {(B.  ,C)} {F_A}$, which implies that $A$ does not cause
either $B$ or $C$, DAG \itref{augsimp1} instead represents the ECI
$\ind {(B,C)} {F_A} A$, expressing the ignorability of the
distribution of $(B,C)$ given $A$.

Note moreover that the Markov equivalence of DAGs \itref{augsimp2} and
\itref{augsimp3} means that we can not interpret the direction of the
arrow between $B$ and $C$ as meaningful in itself.  In particular, in
DAG \itref{augsimp3} the arrow $B\rightarrow C$ does not signify a
causal effect of $B$ on $C$: in this approach, causality is {\em
  only\/} described by ECI properties, and their moralisation-based
graphical representations.
\end{ex}

\begin{ex}
  \label{ex:instraug}
  To endow the instrumental variable problem of \figref{instr} with
  causal content---specifically, relating to the causal effect of $X$
  on $Y$---we might replace it by the augmented DAG of
  \figref{instrf}, where the node $F_X$ now allows us to consider an
  intervention on $X$.

  \begin{figure}[htbp]
  \centering
  \includegraphics [width=.4\textwidth,clip] {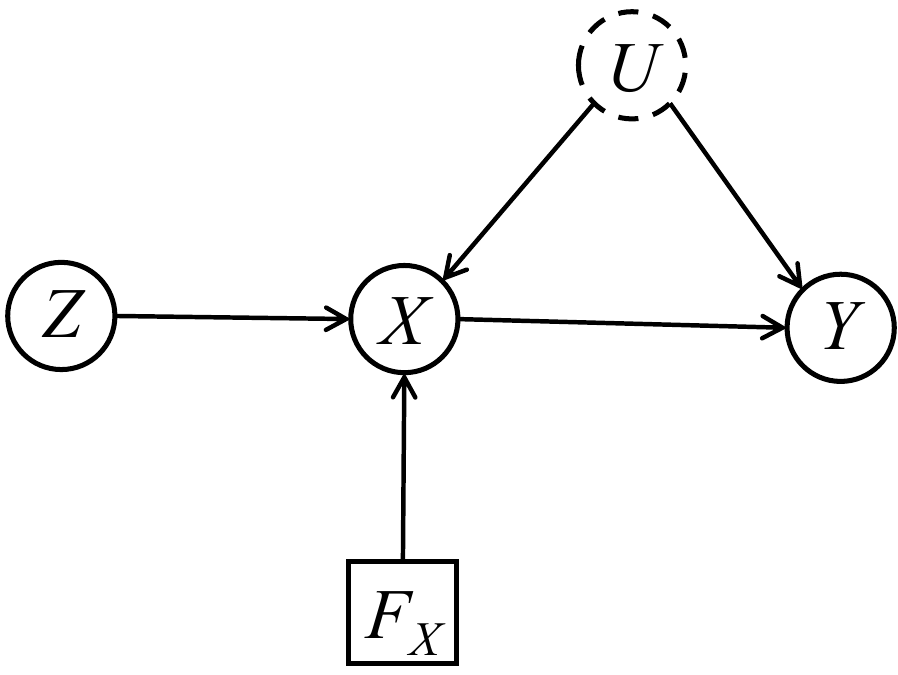}
  \caption{Instrumental variable: Augmented DAG}
  \label{fig:instrf}
\end{figure}

This DAG still represents the probabilistic conditional independencies
\eqref{instrp1} and \eqref{instrp2}, in any regime.  But now it
additionally represents genuine causal properties:
\begin{eqnarray}
  \label{eq:instrc1}
  (U,Z) &\cip& F_X\\
  \label{eq:instrc2}
  Y &\cip& F_X \mid (Z,U,X) 
\end{eqnarray}
Property \eqref{instrc1} says that $X$ has no causal effect on $U$ and
$Z$, these having the same joint distribution in all regimes
(including the idle regime).  Property \eqref{instrc2} entails the
modular conditional ignorability property that the distribution of $Y$
given $(z,u,x)$ (which in fact depends only on $(u,x)$, by
\eqref{instrp1}) is the same, both in the interventional regime where
$X$ is set to $x$, and in the observational regime where $X$ is not
controlled.  Although rarely stated so explicitly, these assumptions
are basic to most understandings of an instrumental variable problem,
and its analysis (which we shall not however pursue here).

If we wanted, we could work with an augmented version of the SPM
representation of \figref{instre}, as in \figref{instref}.
  \begin{figure}[htbp]
  \centering
  \includegraphics [width=.4\textwidth,clip] {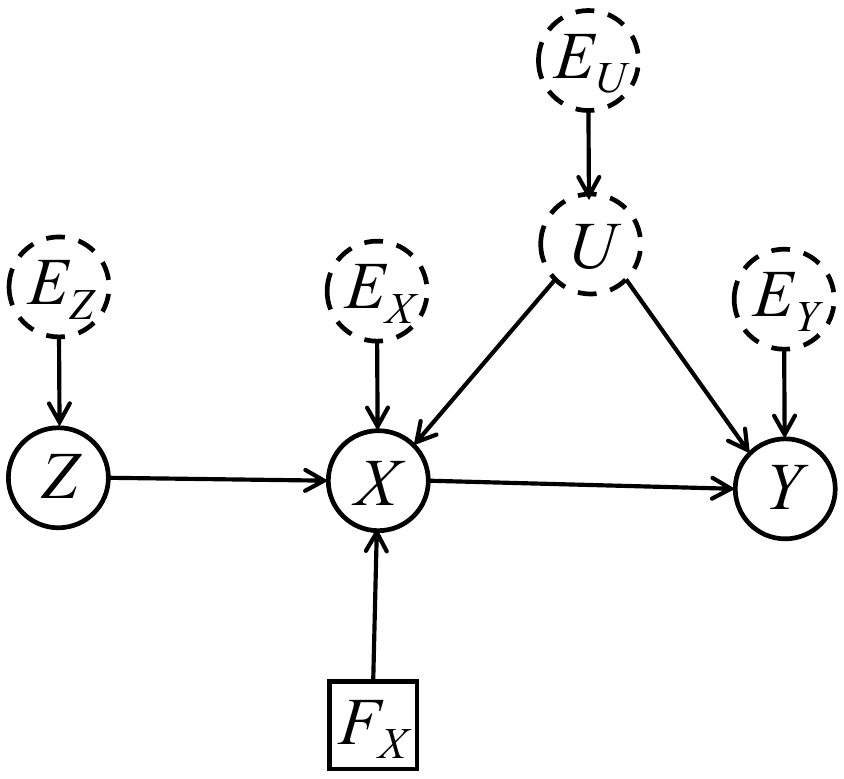}
  \caption{Instrumental variable: Augmented DAG with error variables}
  \label{fig:instref}
\end{figure}
This entails exactly the same ECI properties as \figref{instrf} for
the domain variables and the intervention variable.  With suitably
chosen distributions for the error variables and functional dependence
of each domain variable on its parents, we can recover the same joint
distribution for the domain variables, in any regime, as in the
original probabilistic augmented DAG.  Inclusion of the extra
structure brings nothing new to the table.

Yet another representation of the problem is given in
\figref{instrvf}.
  \begin{figure}[htbp]
  \centering
  \includegraphics [width=.4\textwidth,clip] {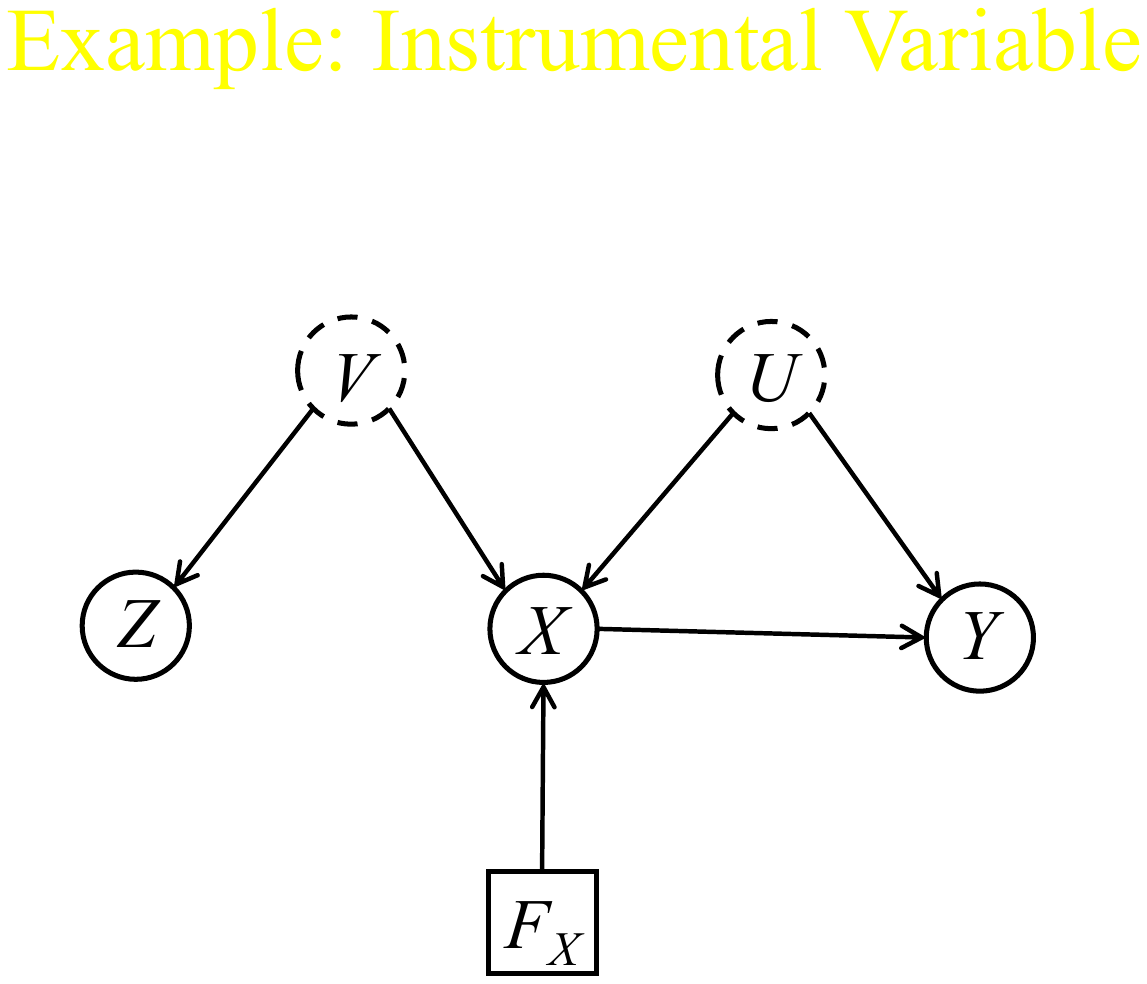}
  \caption{Instrumental variable: Alternative augmented DAG}
  \label{fig:instrvf}
\end{figure}
Here $V$ denotes an additional unobserved variable of no direct
interest.  It can again be checked that both \figref{instrf} and
\figref{instrvf} represent the identical ECI properties between the
variables of interest, $X, Y, Z, U$ and $F_X$.  This identity shows
that the {\bf arrow} $Z\rightarrow X$ in \figref{instrf} should not be
taken as signifying a {\tt direct causal effect} of $Z$ on $X$: we
could equally well regard $X$ and $Z$ as being associated through a
{\tt common cause}, $V$.  \cite{hernan:06} regard \figref{instrf} and
\figref{instrvf} as essentially different---as indeed they would be if
interpreted in these informal terms---and conclude (correctly) that it
is not necessary, for analysing the instrumental variable problem, to
require that $Z$ have a {\tt direct effect} on $X$.  From our point of
view there is no essential difference between \figref{instrf} and
\figref{instrvf}, since even in \figref{instrf} the arrow
$Z\rightarrow X$ should not be interpreted causally.
\end{ex}

\section{Pearlian DAG}
\label{sec:pearl}
\begin{figure}[htbp]
  \centering
  \includegraphics [width=.4\textwidth,clip] {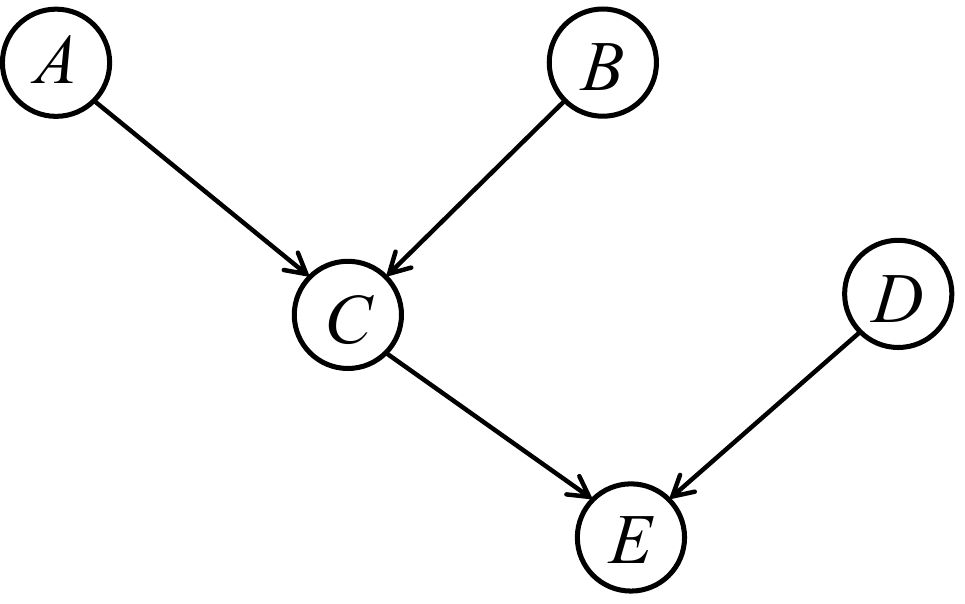}
  \caption{Pearlian DAG}
  \label{fig:pearldag0}
\end{figure}

C  onsider the DAG of \figref{pearldag0}.

As it stands this looks like a probabilitic DAG, representing purely
conditional independence properties, such as $\ind C D {(A,B)}$.
\cite{pearl:book}, however, would endow it with additional causal
meaning, using an interventional interpretation as in \secref{int}.
He would regard it as asserting that, for any node, its conditional
distribution, given its parents, would be the same, in a purely
observational setting, and in any interventional setting that sets the
values of some or all of the variables other than itself.  For
example, it requires:
\begin{quote}  
  \it The distribution of $C$, given $(A,B)$, does not depend on whether
  $A$ and $B$ arose naturally, or were set by external intervention.
\end{quote}

While this is a perfectly clear formal interpretation of
\figref{pearldag0}, it is problematic in that, if we are just
presented with that DAG, we may not know whether it is meant to be
interpreted as representing only probabilistic properties, or is
supposed to be further endowed with Pearl's causal semantics.  This
ambiguity can lead to confusion: a particular danger is to see, or
construct, a probabilistic DAG, and then slide, unthinkingly, into
giving it an unwarranted Pearlian causal interpretation.

\begin{figure}[htbp]
  \centering
  \includegraphics [width=.6\textwidth,clip] {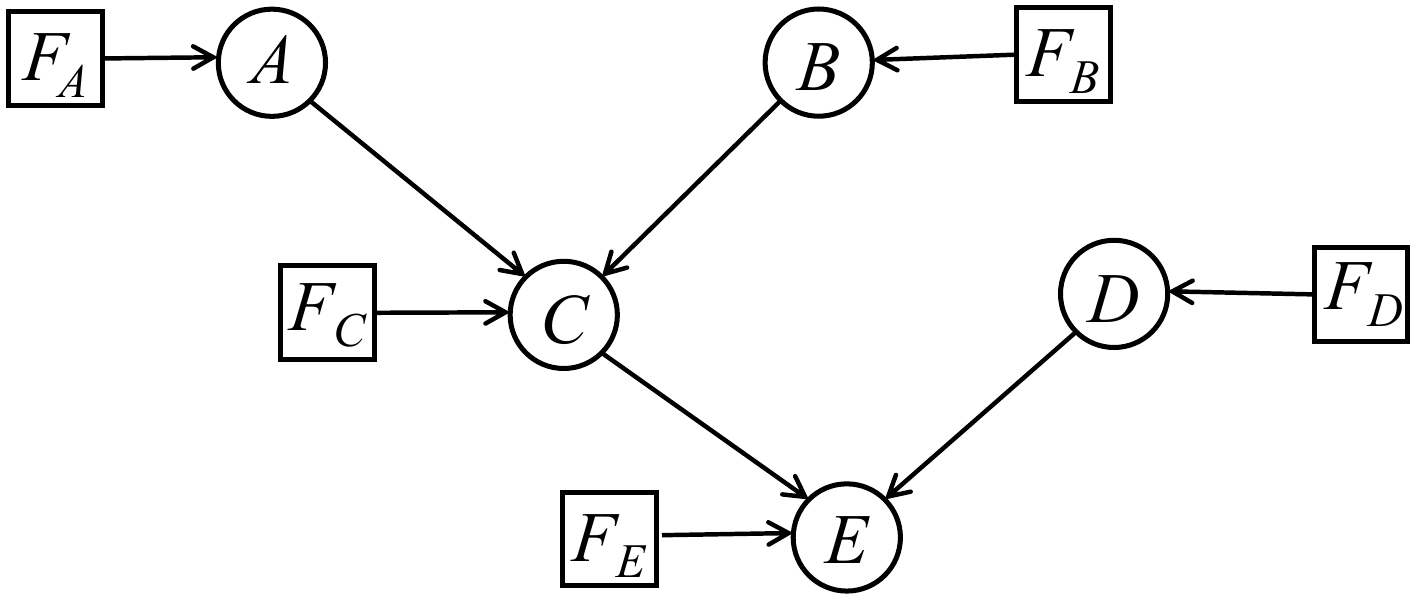}
  \caption{Augmented Pearlian DAG}
  \label{fig:pearldag1}
\end{figure}

We can avoid this problem by explicit inclusion of regime indicators,
one for each domain variable, as for example in \figref{pearldag1}.
Not only is this clearly not intended as a probabilistic DAG, but the
Pearlian causal semantics, which in the case of \figref{pearldag0}
require external specification, are now explicitly represented in the
augmented DAG, by moralisation.  For example, \figref{pearldag1}
represents the ECI $\ind C {(F_A,F_B)} {(A,B,F_C,F_D,F_E)}$.  When
$F_C = c \neq\emptyset$, $C$ has a one-point distribution at $c$, and
this ECI holds trivially.  But for $F_C=\emptyset$ we recover the
property quoted in italics above (under any settings, idle or
interventional, of $F_D$ and $F_E$).

We also note that an augmented Pearlian DAG can have no other Markov
equivalent such DAG, since no arrow can be reversed without creating
or destroying an immorality.  In this sense every arrow now carries
causal meaning.

However, just as we should not automatically interpret a regular DAG
as Pearlian, so we should not unthinkingly simply augment it by adding
an intervention indicator for each domain variable---which would have
the same effect.  We must consider carefully whether the many very
strong properties embodied in any Pearlian or augmented DAG, relating
probabilistic properties (parent-child distributions) across distinct
regimes, are justifiable in the specific applied context we are
modelling.

\subsection{Structural Causal Model}
\label{sec:structcaus}
We can also reintepret a SPM, such as in \figref{instre}, using
Pearlian semantics\footnote{where the possibility of intervention is
  envisaged for each of the domain variables, but not for the
  fictitious error variables}, as a causal model---a {\em Structural
  Causal Model (SCM)}.  This would then assert that the distributions
of the error variables, and the functional dependence of each
unintervened domain variable on its parents, are the same in all
regimes, whether idle or subject to arbitrary interventions.

Again, to avoid confusion with a SPM, it is advisable to display a SCM
as an augmented DAG, by explicitly including intervention indicators
(as in \figref{instref}, but having an intervention indicator
associated with {\em every\/} domain variable.)

However, construction and inclusion of fictitious error variables is
of no consequence, since, if we concentrate on the domain variables
alone, their probabilistic structure, in any regime, will be exactly
the same as for the fully probabilistic augmented DAG.  In particular,
no observations on the domain variables, under any regime or regimes,
can resolve the arbitrariness in the specifications of the error
distributions and the functional relationships in a SCM.  On this
second rung of the ladder of causation, the additional structure of
the SCM once again gains us nothing.

\section{Causes of Effects}
\label{sec:coe}
So far we have only considered modelling and analysing problems
relating to the ``effects of causes (EoC)'', where we consider the
downstream consquenced of an intervention.  An entirely different
class of causal questions relates to ``causes of effects (CoE)''
\citep{apd/mm:eoccoe}, where, in an individual case, both a putative
causal variable $X$ and a possible consequence $Y$ of it have both been
observed, and we want to investigate whether $X$ was indeed the
cause of $Y$.  Such problems arise in the legal context,
for example when an individual sues a pharmaceutical company claiming
that it was because she took the company's drug that she developed a
harmful condition.

At the simplest level we may only have information about the joint
distribution of $(X,Y)$ and their values, $x$ and $y$, for the case at
hand.  But this is not enough to address the CoE question, which
refers, not to an unknown fact or variables, but to an unknown
relationship: was it causal?  To try to undersand this question takes
us into new territory---the third rung of the ladder of causation.

Although by no means totally satisfactory, this question is most
commonly understood as relating to a counterfactual situation.
Suppose both $X$ and $Y$ are binary, and we have observed $X=Y=1$.  We
mught express the ``event of causation'' as
\begin{quote}
  \em The outcome variable would have been different (\ie, $Y=0$) if
  the causal variable had been different (\ie, $X=0$)
\end{quote}
But the hypothesis here, $X=0$, contradicts the known fact that
$X=1$---it is {\em counterfactual\/}.  Since in no case can we observe
both $X=1$ and $X=0$---what has been called ``the fundamental problem
of causal inference'' \citep{pwh:jasa}---it would seem impossible
to address this question, at any rate, on the basis of purely factual
knowledge of the joint distribution of $(X,Y)$.  So a more complicated
framework is required, necessitating more complex assumptions and
analyses \citep{apd:kent}.

One approach builds on the idea of ``potential responses'',
popularised by \cite{dbr:jep,dbr:as} initally for addressing EoC
questions---which, as we have seen, can progress perfectly well
without them.  They do however seem essential for formulating
counterfactual questions.  For the simple example above, we duplicate
the response $Y$, replacing it by the pair $(Y_0,Y_1)$, with $Y_x$
conceived of as a potential response, that would be realised if in
fact $X=x$.  Then the {\em probability of causation\/}, \pc, can be
defined as the conditional probability, given the data, of the event
of causation:
\begin{equation}
  \label{eq:pc}
  \pc = \pr(Y_0=0\mid X=1, Y_1=1).
\end{equation}
There is however a difficulty: on account of the fundamental problem
of causal inference, no data, of any kind, could ever identify the
joint distribution of $(Y_0,Y_1)$, so $\pc$ is not estimable.  It
turns out that data supplying the distribution of the observable
variables $(X,Y)$ can be used to set interval bounds on \pc, and these
bounds can sometimes be refined if we can collect data on additional
variables \citep{apd/mm/rm,apd/mm:sef}, but only in very special cases
can we obtain a point value for \pc.

\subsection{Graphical representation}
\label{sec:graph-rep}
Because \pc\ cannot be defined in terms of observable domain
variables, we cannot represent a CoE problem by means of a regular or
augmented DAG on these variables.  It is here that the expanded
SCM version appears to come into its own.

\begin{figure}[htbp]
  \centering
  \includegraphics [width=.4\textwidth,clip]{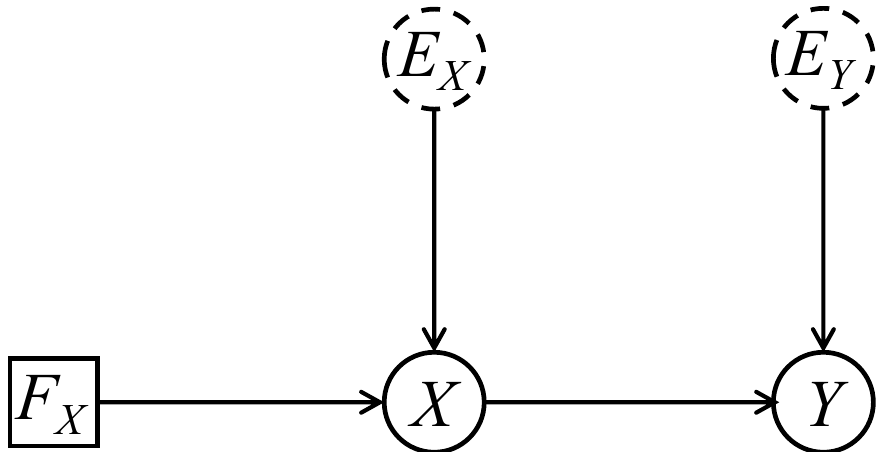}
  \caption{Simple SCM}
  \label{fig:simscm}
\end{figure}

Thus consider the  simple SCM of \figref{simscm}.
In it, we have $X=f_X(E_X)$, $Y=f_Y(X,E_Y)$, with $\indo {E_Y} X$.
The deterministic structure of a SCM supports the introduction of
potential responses (something that can not be done with a purely
probabilistic DAG): $Y_x = f_Y(x,E_Y)$.\footnote{Indeed, we could
  replace $E_Y$ with the pair $\bY := (Y_0,Y_1)$, and the function
  $F_Y$ with the ``look-up'' function $L_Y$, where
  $L_Y(x, \bY) = Y_x$.}  Then $X, Y_0, Y_1$ are all functions of
$(E_X, E_Y)$, and thus have a joint distribution entirely determined
by that of $(E_X, E_Y)$.  And given this joint distribution we can
compute \pc in \eqref{pc}.  Likewise, given a more general SCM we can
compute answers to CoE questions.

It would seem then that this approach bypasses the difficulties
alluded to above.  This, however, is illusory, since such a solution
is available only when we have access to a fully specified SCM.  As
discussed in \secref{struct} and \secref{structcaus}, there will be
many distinct PSMs or SCMs consistent with a given probabilistic or
augmented DAG model for the domain variables.  Since the probabilistic
structure is the most that we can learn from empirical data, we will
never be able to choose between these distinct SCM versions of the
problem.  However, different SCMs will lead to different answers to
the CoE questions we put to them.  Thus for \figref{simscm}, we have
$\pc = \pr(Y_0=0\mid Y_1=1)$, since $\indo \bY X$, and this will
depend on the dependence between $Y_0$ and $Y_1$ as embodied in the
SCM.  But because of the fundamental problem of causal inference, this
dependence can never be identified from empirical data.  So different
SCMs inducing the same probabilistic structure, which are entirely
indistinguishable empirically, will lead to different answers.  When
we allow for all possible such SCMs, we are led back to the interval
bounds for $\pc$ discussed above.

\section{Discussion}
\label{sec:disc}
We have surveyed a variety of directed graphical models, with varying
syntax (including or omitting error variables or regime indicators)
and semantics (formal or informal, modelling probabilistic or causal
situations).  Different semantics are relevant to different rungs of
the ladder of causation.

When presented with a simple DAG it may not be obvious how it is
supposed to be interpreted, and there is an ever-present danger of
misinterpretation, or of slipping too easily from one interpretation
(\eg, purely probabilistic) to another (\eg, causal).  This can
largely be avoided by always using a simple DAG to model a
probabilistic problem (on the first rung of the ladder), and an
augmented DAG to model a causal problem (on the second rung).  In both
cases, the moralisation procedure provides the semantics whereby
interpretive properties can be read off the graph.

DAGs such as SCMs, that involve, explicity or implicitly, error
variables and functional relationships, can be used on all three rungs
of the ladder.  However they can not be identified empirically.  For
rungs one and two this is unimportant, and all equivalent versions,
inducing the same underlying purely probabilistic DAG, will yield the
same answers as obtainable from that underlying DAG.  For rung 3,
which addresses the probability of causation in an individual
instance, only an approach based on SCMs is available. However,
different but empirically indistinguishable SCMs now deliver different
answers to the same causal question.  Taking this into account we may
have to be satisfied with an interval bound on the desired, but
unidentifiable, probability of causation.

\bibliographystyle{apalike}
\bibliography{strings,causal,ci,allclean}

\end{document}